\documentclass[twoside,10pt,leqno]{amsart}
\usepackage{amsfonts}
\usepackage{amsmath}
\usepackage{amscd}
\usepackage{amssymb}
\usepackage{amsthm}
\usepackage{amsrefs}
\usepackage{latexsym}
\usepackage{bbm}
\numberwithin{equation}{section}
\begin{document}
\newtheorem{theorem}{Theorem}
\newtheorem{lemma}{Lemma}
\newtheorem{prop}{Proposition}
\newtheorem{corollary}{Corollary}
\newtheorem{conjecture}{Conjecture}
\numberwithin{equation}{section}
\newcommand{\dif}{\mathrm{d}}
\newcommand{\intz}{\mathbb{Z}}
\newcommand{\ratq}{\mathbb{Q}}
\newcommand{\natn}{\mathbb{N}}
\newcommand{\comc}{\mathbb{C}}
\newcommand{\rear}{\mathbb{R}}
\newcommand{\prip}{\mathbb{P}}
\newcommand{\uph}{\mathbb{H}}

\title{An Improvement for the Large Sieve for Square Moduli}
\author{Stephan Baier \and Liangyi Zhao}
\date{\today}
\maketitle

\begin{abstract}
We establish a result on the large sieve with square moduli.  These bounds improve recent results by S. Baier \cite{Ba1} and L. Zhao \cite{Zha}.
\end{abstract}

\noindent {\bf Mathematics Subject Classification (2000)}: 11B57, 11L07, 11L40, 11N35, 11L15. \newline

\noindent {\bf Keywords}: large sieve, Farey fractions in short intervals, estimates on exponential sums and integrals.

\section{Introduction and Statement of the Main Results}

Large sieve was an idea originated by J. V. Linnik \cite{JVL1} in 1941 while studying the distribution of quadratic non-residues.  Refinements of this idea were made by many.  In this paper, we develop an improvement for large sieve inequality for square moduli.  More in particular, we aim to have an estimate for the following sum
\begin{equation} \label{A1}
\sum\limits_{q\le Q} \sum\limits_{\substack{ a=1\\(a,q)=1}}^{q^2}
\left| \sum_{n= M+1}^{M+N} a_n e\left(\frac{a}{q^2}n\right)
\right|^2,
\end{equation}
where henceforth $M \in \intz$, $Q, N \in \natn$ and $\{ a_n \}$ is an arbitrary sequence of complex numbers.  In the sequel, we set
\[ S(x): =  \sum_{n=M+1}^{M+N} a_n e\left(xn\right), \; \mbox{and} \; Z:=\sum\limits_{n=M+1}^{M+N} \vert a_n\vert^2. \]

With $q^2$ replaced by $q$ in \eqref{A1}, it is
\begin{equation} \label{kis1}
 \ll (Q^2+N)Z.
\end{equation}
This is in fact the consequence of a more general result first introduced by H. Davenport and H. Halberstam \cite{DH1} in which the Farey fractions in the outer sums of \eqref{A1} can be replaced by any set of well-spaced points.  Montgomery and Vaughan \cite{MVa} showed that the $\ll$ can be replaced by $\leq$ in \eqref{kis1}.  Literature abound on the subject of the classical large sieve. See \cites{BD1, HD, DH1, PXG, JVL1, Lem, HM2, MVa}.   Applying the said more general result, \eqref{A1} is bounded above by 
\begin{equation}
\ll (Q^3+QN)Z, \mbox{ and } \ll (Q^4+N)Z\label{A4}
\end{equation}
(See \cites{SBLZ, Zha}).  In \cite{Zha} it was proved that the sum \eqref{A1} can be estimated by 
\begin{equation}
\ll \log (2Q) \left(Q^3+\left(N\sqrt{Q}+ \sqrt{N}Q^2\right)N^{\varepsilon}\right)Z,\label{A5}
\end{equation}
where the implied constant depends on $\varepsilon$.  Also in \cite{Ba1}, it was shown that \eqref{A1} is
\begin{equation} \label{Baiertheorem}
\ll (QN)^{\varepsilon} \left( Q^3+ N + \sqrt{N}Q^2 \right)Z.
\end{equation}
Moreover, in the same paper \cite{Ba1}, Baier showed that \eqref{A1} is 
\[ \ll \left\{ \begin{array}{cc} Q^{3/5+\varepsilon}NZ, & \mbox{if} \; Q \leq N^{5/12}, \\ Q^{3+\varepsilon}Z, & \mbox{if} \; Q > N^{5/12}. \end{array} \right. \]
It was conjectured in \cite{Zha} that \eqref{A1} is
\begin{equation} \label{conj}
\ll Q^{\varepsilon} (Q^3+N)Z.
\end{equation}
Furthermore, the inequality for higher power moduli was also studied in \cites{Zha, SBLZ}. \newline

In this paper, we prove the following.

\begin{theorem} \label{squaremod1}
With $\varepsilon >0$ given and under the usual notations, we have
\begin{equation} \label{squaremod1eq}
\sum\limits_{q\le Q} \sum\limits_{\substack{a=1 \\ (a,q)=1}}^{q^2}
\left\vert S\left(\frac{a}{q^2}\right)\right\vert^2 \ll (QN)^{\varepsilon}
\left(Q^3+N+\min\{N\sqrt{Q},\sqrt{N}Q^2\}\right)Z,
\end{equation}
where the implied constant depends on $\varepsilon$ alone.
\end{theorem}

This theorem provides a better majorant in the range $N^{1/3+\varepsilon} \leq Q \leq N^{5/12-\varepsilon}$ and establishes the conjecture in \eqref{conj} for the range $Q \gg N^{2/5}$.  Theorem~\ref{squaremod1} reduces to counting Farey fractions with square denominators in short intervals.  We do this counting in two ways.  The first is a generalization of the techniques in \cite{Zha} using estimates for certain Weyl sums, and the second is a more refined treatment of certain exponential integrals in \cite{Ba1}.  The improvement comes from using different estimates for the integrals depending on the location of the stationary points of the integrals of interest and a more explicit evaluation of the quadratic Gauss sums from \cite{Ba1}.  Theorem~\ref{squaremod1} follows when the two estimates are combined and compared with previously known results. \newline

Furthermore, as an application of the large sieve for square moduli, we note the following Bombieri-Vinogradov type bound for prime squares.  We have that, using \eqref{Baiertheorem} and arguments similar to those in the proof of the
classical Bombieri-Vinogradov theorem, 
\[ \sum_{p \leq x^{2/9-\varepsilon}} p \max_{y\leq x} \max_{\substack{a \bmod{p^2} \\ p \nmid a}} \left| \psi(y; p^2,a) - \frac{y}{\varphi(p^2)} \right| \ll \frac{x}{(\log x)^A}, \]
for any $A>0$.  The restriction to only prime square moduli provides
additional convenience in the reduction to primitive characters and the
authors also believe that it is worthwhile to investigate results of
this type for general sparse sets of moduli. \newline

We note here that a Bombieri-Vinogradov type theorem more general than the above in which all square moduli are considered, not just the prime square moduli, has appeared in \cite{SBLZ3}
while the present paper was under review.  Moreover, the more general result improves
some earlier results of Elliott \cite{PDTAE} and Mikawa and Peneva \cite{HMTPP}.

\section{Preliminary Lemmas}

In this section, we quote some lemmas that we shall use later.  We shall need the following general version of the large sieve inequality.
\begin{lemma} \label{generalls}
Let $\left\{ \alpha_r\right\}_{r=1}^R$ be a
sequence of real numbers. Suppose that $0<\Delta\le 1/2$ and 
$R\in \mathbbm{N}$. Set 
$$
K(\Delta):=\max\limits_{\alpha\in \mathbbm{R}} 
\sum\limits_{\substack{r=1 \\
\| \alpha_r -\alpha \| \le \Delta }}^R 1,
$$
where $\| x \|$ denotes the distance of a real $x$ to its closest integer.
Then 
$$
\sum\limits_{r=1}^R \left\vert S\left(\alpha_r\right)\right\vert^2
\ll K(\Delta)(N+\Delta^{-1})Z,
$$
with an absolute $\ll$-constant.
\end{lemma}
\begin{proof}
This is Theorem 2.11 in \cite{Lem}.
\end{proof}

We also need the Poisson summation formula.

\begin{lemma} [Poisson Summation Formula] \label{poisum}
Let $f(x)$ be a function on the real numbers that is piece-wise continuous with only finitely many discontinuities and for all real numbers, $a \in \rear$, satisfies
\begin{equation*}
f(a) = \frac{1}{2} \left[ \lim_{x \to a-} f(x) + \lim_{x \to a+} f(x) \right].
\end{equation*}
Moreover, $f(x) \ll (1+|x|)^{-c}$ for some $c>1$ with an absolute implied constant.  Then we have
\begin{equation*}
\sum_{n=-\infty}^{\infty} f(n) = \sum_{n=-\infty}^{\infty} \hat{f}(n), \; \mbox{where} \; \hat{f}(x) = \int_{-\infty}^{\infty} f(y) e(xy) \dif y,
\end{equation*}
the Fourier transform of $f(x)$.
\end{lemma}
\begin{proof} This is quoted in \cite{Bum}. \end{proof}

We shall not succeed in proving our contention without the following lemma.

\begin{lemma}[Weyl Shift] \label{weylk}
Let $I$ be an interval of length $N$ and $f(x)$ be a polynomial of degree $k \geq 2$ with real coefficients.  Set $\kappa=2^{k-1}$ and let the leading coefficient, the coefficient of $x^k$, of $f(x)$ be $\alpha$.  Also set
\begin{equation*}
S= \sum_{n \in I} e(f(n)).
\end{equation*}
Then we have
\begin{equation*}
|S|^{\kappa} \leq 2^{2\kappa}
N^{\kappa-1} + 2^{\kappa} N^{\kappa-k} \sum_{r_1,\cdots,r_{k-1}} \min \left( N, \frac{1}{\| \alpha k! r_1 \cdots r_{k-1} \|} \right),
\end{equation*}
where each $r$ runs from 1 to $N-1$ and $\| x \| = \min_{l\in \intz} |x-l|$.
\end{lemma}
\begin{proof} This is Lemma 5.6 in \cite{ET}. \end{proof}

We also need the following asymptotic formula for exponential integrals with weights.

\begin{lemma}[Stationary Phase with Weights] \label{Ivictheorem}
Let $f(z)$, $g(z)$ be two functions of the complex variable $z$ and
$[a,b]$ a real interval such that:
\begin{enumerate}
\item For $a\le x\le b$ the function $f(x)$ is real and $f''(x)>0$.
\item For a certain positive differentiable function $\mu(x)$, defined on
$a\le x\le b$, $f(z)$ and $g(z)$ are analytic for $a\le x\le b$,
$|z-x| \le \mu(x)$.  
\item There exist positive functions $F(x)$, $G(x)$ defined on $[a,b]$ such 
that for $a\le x\le b$, $|z-x|\le \mu(x)$ we have
$$
g(z)\ll G(x),\ \ \ \ \ \ f'(z)\ll \frac{F(x)}{\mu(x)},\ \ \ \ \ \ 
|f''(z)|^{-1}\ll \frac{\mu(x)^2}{F(x)},
$$
and the $\ll$-constants are absolute.
\end{enumerate}

Let $k$ be any real number, and if $f'(x)+k$ has a zero in [a,b] denote it by
$x_0$. Let the values of $f(x)$, $g(x)$, and so on, at $a$, $x_0$, and $b$ be
characterized by the suffixes $a$, $0$, and $b$ respectively. Then, for some
absolute constant $C>0$, we have
\begin{equation*}
\begin{split}
\int_a^b g(x) e (f(x)+kx) \dif x = \frac{g_0}{\sqrt{f_0''}}&e \left( f_0+kx_0+\frac{1}{8} \right) + O \left(\int\limits_a^b G(x)\exp\left(-C|k|\mu(x)-CF(x)\right)
({\rm d}x+|{\rm d}\mu(x)|)\right) \\
& + O \left(G_0\mu_0F_0^{-3/2} + G_a\left(|f_a'+k|+\sqrt{f_a''}\right)^{-1}+ G_b\left(|f_b'+k|+\sqrt{f_b''}\right)^{-1}\right).
\end{split}
\end{equation*}
\end{lemma}
\begin{proof}
This is Theorem 2.2. in \cite{Ivic}.
\end{proof}

Then we need the following lemma for the estimation of exponential integrals.

\begin{lemma} \label{gkest}
Assume that $f$ and $g$ are real-valued, twice differentiable functions on $[a,b]$.  Also assume that $g/f'$ is monotonic and that 
\[ \left| \frac{f'(x)}{g(x)} \right| \geq \lambda. \]
Then we have
\[ \int_a^b g(x) e(f(x)) \dif x \ll \lambda^{-1}. \]
\end{lemma}
\begin{proof}
This is Lemma 3.1 in \cite{Kol}.
\end{proof}

We shall also need to transform and estimate certain Gauss sums.  We define
\begin{equation} \label{gausssumdef}
G (a,l;c) = \sum_{d \pmod{c}} e \left( \frac{ad^2+ld}{c} \right).
\end{equation}
We have the following lemma.

\begin{lemma} \label{gausssumtrans}
Assuming that $(a,c)=1$ and $a\overline{a}\equiv 1$ mod $c$, we have
\begin{enumerate}

\item If $l$ is even, then
\[ G(a,l;c) = e \left( - \frac{\overline{a}l^2}{4c} \right) G(a,0;c); \]

\item If $l$ is odd, then
\[ G(a,l;c) = e \left( - \frac{\overline{a}(l^2-1)}{4c} \right) G(a,1;c). \]

\end{enumerate}
\end{lemma}

\begin{proof}
This is Lemma 7.11 in \cite{Kol}.
\end{proof}

To estimate these Gauss sums, we use the following Lemma.

\begin{lemma} \label{gausssumest}
If $(a,c)=1$, then we have
\[ |G(a,l;c)| \leq 2\sqrt{c}. \]
\end{lemma}
\begin{proof}
This is (7.4.2) in \cite{Kol}.
\end{proof}

\section{Counting Farey fractions in short intervals}

We aim to estimate
\begin{equation} \label{aim}
\sum_{q \leq Q} \sum_{\substack{a=1 \\ (a,q)=1}}^{q^2} \left| S \left( \frac{a}{q^2} \right) \right|^2. 
\end{equation}
It suffices to consider the above expression with the outer-most sum 
running over dyadic intervals.  To that end, we begin with the following modified version of Lemma~\ref{generalls}.

\begin{lemma} \label{modified}
Let $\left\{ \alpha_r\right\}_{r=1}^R$ and $\left\{ \beta_l \right\}_{l=1}^L$ be two sequences of real numbers.  Suppose that $0<\Delta\le 1/2$ and for every $\alpha \in \rear$ there exists $\beta_l$ with $1 \leq l \leq L$ such that
\[ \| \beta_l - \alpha \| \leq \Delta.\]
Put 
$$
K'(\Delta):=\max_{1 \leq l \leq L} 
\sum_{\substack{r=1 \\ \| \alpha_r -\beta_l \| \le \Delta}}^R 1.
$$
Then 
$$
\sum\limits_{r=1}^R \left| S\left(\alpha_r\right)\right|^2
\ll K'(\Delta)(N+\Delta^{-1})Z,
$$
with an absolute $\ll$-constant.
\end{lemma}
\begin{proof}
This follows easily from Lemma~\ref{generalls}.
\end{proof}

In our situation, let $\alpha_1,...,\alpha_R$ be the sequence of Farey fractions $a/q^2$ with $Q < q \leq 2Q$, $1\le a\le q$ and $(a,q)=1$.  We shall use Lemma~\ref{modified} to estimate 
\[ \sum\limits_{r=1}^R \left| S\left(\alpha_r\right)\right|^2 \]
and choose the $\beta_l$'s in an appropriate way.  First we set
\begin{equation}
\tau:=\frac{1}{\sqrt{\Delta}}.\label{P1}
\end{equation}
Let $\beta_1, \ldots, \beta_L$ be 
\begin{equation} \label{betal}
 \frac{b}{r} + \frac{1}{kr^2}
\end{equation}
with $r \in \natn$ and $r \leq \tau$, $(b,r)=1$, $1 \leq b \leq r$ and $k \in \intz$ with $\lceil r^{-1} \Delta^{-1/2} \rceil\leq |k| \leq 
\lceil r^{-2} \Delta^{-1} \rceil$.  Here and after, we set $\lceil x \rceil = \min_{l \in \intz} \{ l \geq x \}$ for $x \in \rear$.  We want to show that the $\beta_l$'s above satisfy the conditions of Lemma~\ref{modified}.  By Dirichlet's approximation theorem, every $\alpha \in \rear$ can be written in the form
\begin{equation} \label{P2}
\alpha=\frac{b}{r}+z, \ \ \mbox{ where }\ \  r\le \tau,\ (b,r)=1,\ 
\vert z\vert \le \frac{1}{r\tau}.
\end{equation}
We must show that for every $|z| \leq (r\tau)^{-1}$, there is a $k \in \intz$ with $\lceil r^{-1} \Delta^{-1/2} \rceil
\leq |k| \leq \lceil r^{-2} \Delta^{-1}\rceil$ such that
\[ \left| z - \frac{1}{kr^2} \right| \leq \Delta. \]
First, we have if $0 \leq z \leq (\kappa r^2)^{-1}$ with $\kappa= \lceil r^{-2}\Delta^{-1}\rceil$, then 
\[ \left| z - \frac{1}{\kappa r^2} \right| \leq \Delta. \]
Otherwise, if $z>(\kappa r^2)^{-1}$, it suffices to show that
\begin{equation} \label{ineq1}
\frac{1}{kr^2} - \frac{1}{(k+1)r^2} \leq \Delta,
\end{equation}
for the $k$'s in question and for $K= \lceil r^{-1} \Delta^{-1/2} \rceil$,
\begin{equation} \label{ineq2}
 \left| \frac{1}{r\tau} - \frac{1}{Kr^2} \right| \leq \Delta.
\end{equation}
The left-hand side of \eqref{ineq1} is 
\[ \frac{1}{k(k+1)r^2} \leq \frac{1}{k^2r^2} \leq \Delta, \]
provided $k \geq K$.  Furthermore, we have
\[ \frac{1}{((r\sqrt{\Delta})^{-1}+1)r^2} \leq \frac{1}{Kr^2} \leq \frac{1}{r\tau} \]
and thus the left-hand of \eqref{ineq2} is
\[ \leq \left| \frac{1}{r\tau} - \frac{1}{((r\sqrt{\Delta})^{-1}+1)r^2} \right| = \frac{\Delta}{1+r\sqrt{\Delta}} \leq \Delta. \]
For $z<0$, the arguments are similar.  \newline

For $\alpha\in \mathbbm{R}$ we put
$$
P(\alpha):= \sum_{\substack{Q< q \leq 2Q, (a,q)=1 \\
| a/q^2 -\alpha | \leq \Delta}} 1.
$$
Then we have
\begin{equation} \label{60}
K'(\Delta) \leq 2\max\limits_{1 \leq l \leq L} P(\beta_l).
\end{equation}

Summarizing the above observations, we deduce the following.

\begin{lemma} \label{engels}
We have  
\begin{equation}
K'(\Delta)\le 2\max_{\substack{ r\in \mathbbm{N} \\ r\le \tau }} 
\max_{\substack{ b\in \mathbbm{Z} \\ (b,r)=1 }} \max_{K \leq |k| \leq \kappa} 
 P\left(\frac{b}{r}+ \frac{1}{kr^2} \right).
\end{equation}
\end{lemma}

Therefore, by the virtue of the preceding lemma, it suffices to estimate $P\left( \alpha \right)$ for $\alpha$ with 
\begin{equation} \label{alphacond}
\alpha = b/r +z,\ \ \ \ \ \ r\leq \tau,\ (b,r)=1,\ z=1/(kr^2),
\end{equation} 
where $k \in \intz$ and $\lceil r^{-1} \Delta^{-1/2} \rceil \leq |k| 
\leq \lceil r^{-2} \Delta^{-1}\rceil$.  We note that $\alpha$ satisfies (\ref{P2}) if it satisfies (\ref{alphacond}). Moreover, it is enough to consider only $k>0$ since
\[ P \left( \frac{b}{r} + \frac{1}{kr^2} \right) = P \left( \frac{r-b}{r} - \frac{1}{kr^2} \right).\]
Consequently, we assume henceforth that
\begin{equation} \label{otto}
z=1/(kr^2) \; \mbox{with} \; k \in \natn \; \mbox{and} \; 
\lceil r^{-1} \Delta^{-1/2} \rceil \leq k \leq \lceil r^{-2} \Delta^{-1}
\rceil.
\end{equation}

\section{A first estimate for $P(\alpha)$}

First, we henceforth set $\phi (x) = \left( \frac{\sin \pi x}{2x} \right)^2$, a constant multiple of F\'ejer kernel.  We note that $\phi (x)$ is non-negative, $\phi(x) \geq 1$ for $|x| \leq \frac{1}{2}$ and $\phi (0) = \pi^2/4$.  Therefore
\begin{equation} \label{estm}
P ( \alpha ) \leq \sum_{Q < q \leq 2Q} \sum_{a=-\infty}^{\infty} \phi \left( \frac{a-\alpha q^2}{8Q^2\Delta} \right).
\end{equation}
We find it most convenient to choose $\phi(x)$ this way, since its Fourier transform is a function of compact support, specifically 
\[ \hat{\phi}(s)=\frac{\pi^2}{4} \max(1-|s|,0). \]

Now applying Poisson summation, Lemma~\ref{poisum}, with a linear change of variable, to \eqref{estm}, we get that it is
\begin{equation}
8Q^2 \Delta \sum_{Q < q \leq 2Q} \sum_{j=-\infty}^{\infty} \hat{\phi} (8Q^2\Delta j) e \left( \alpha j q^2 \right) .
\end{equation}

More precisely, the above is
\[ 8Q^2 \Delta \sum_{|j| < (8Q^2\Delta)^{-1}} \sum_{Q < q \leq 2Q} \left( 1 - 8|j| Q^2\Delta \right) e \left( \alpha j q^2 \right) \ll  Q^3 \Delta + Q^2\Delta \sum_{0< j < (8Q^2\Delta)^{-1}} \left| \sum_{Q < q \leq 2Q}  e \left( \alpha j q^2 \right) \right|, \]
where the first term above corresponds to the contribution of $j=0$.  Applying Cauchy-Schwarz inequality, we see that the square of the above expression is bounded by
\begin{equation*}
\ll Q^6 \Delta^2 + Q^2 \Delta \sum_{0 < j < (8Q^2\Delta)^{-1}} \left| \sum_{Q < q \leq 2Q} e \left( \alpha j q^2 \right) \right|^2.
\end{equation*}
Applying Weyl Shift, Lemma~\ref{weylk}, to the inner-most sum of the second term of the above, we see that the double sum of the said term is
\begin{equation} \label{doh}
\ll \sum_j Q + \sum_j \sum_{0 < l < Q} \min \left\{ Q, \left\| 2 \alpha jl \right\|^{-1} \right\} \ll (Q \Delta)^{-1} + \sum_{0 < m < (Q \Delta)^{-1}} \tau (m) \min \left\{ Q, \left\| \alpha m \right\|^{-1} \right\},
\end{equation}
where $\tau(m)$ is the divisor function, is $O(m^{\varepsilon})$ and estimates the multiplicity of representations of $m=2jl$. \newline

Now we have $\alpha = b/r + z$ with $b, r$ and $z$ as in \eqref{P2}.  The second term in \eqref{doh} is
\begin{equation*}
\begin{split}
& \ll (Q\Delta^{-1})^{\varepsilon} \sum_{0 \leq n < (Q \Delta r)^{-1}} \sum_{m\in (nr, (n+1)r] } \min \left\{ Q, \; \left\| \left( \frac{b}{r} +z \right) m \right\|^{-1} \right\} \\
& = (Q\Delta^{-1})^{\varepsilon} \sum_{0 \leq n < (Q \Delta r)^{-1}} \sum_{l=1}^r \min \left\{ Q, \; \left\| \left( \frac{b}{r} +z \right) (nr+l) \right\|^{-1} \right\}.
\end{split}
\end{equation*}
The inner sum of the above is
\begin{equation*}
\sum_{l=1}^r \min \left\{ Q, \; \left\| znr + l \frac{b}{r} + lz \right\|^{-1} \right\} \ll Q+\sum_{l=1}^{r-1} \frac{r}{l} \ll Q + r \log r.
\end{equation*}
The first of the above inequalities arrives since for each $p/r$ with $0 \leq p \leq r-1$, we find at most three $l$'s such that
\[ \left\| \frac{p}{r} - (znr + l \frac{b}{r} + lz) \right\| \leq \frac{1}{2r}, \]
as
\[ lz \leq rz \leq \sqrt{\Delta} \leq \frac{1}{r}, \]
where we have used \eqref{P2}. \newline

Combining everything, we have the following.

\begin{theorem}
Let $\varepsilon >0$ be given and $\alpha$ satisfy \eqref{P2}.  We have
\begin{equation} \label{firstestPalpha}
P(\alpha) \ll Q^3 \Delta + (Q\Delta^{-1})^{\varepsilon} \left( \frac{Q}{\sqrt{r}} + \sqrt{Q} + Q\sqrt{r\Delta}+Q^{3/2}\sqrt{\Delta}\right),
\end{equation}
where the implied constant depends on $\varepsilon$ alone.
\end{theorem}

\section{Transformation of $P(\alpha)$}
We first note that $z\ge \Delta/2$ by (3.9) and $r\le 1/\sqrt{\Delta}$. If 
$\Delta/2\le z\le \Delta$, then we have
\begin{equation} \label{verysmalldelta}
P(\alpha) \ll \Delta^{-\varepsilon}\left(1+Q_0\Delta r+Q_0^{3/2}\Delta\right)
\end{equation}
by the first inequality of Lemma 6 in \cite{Ba1}. We note that this lemma remains valid, with a different $\ll$-constant, if the condition $z\ge\Delta$ in this lemma is replaced by $z\ge \Delta/2$. 
  
Throughout the following, we assume that $z\ge \Delta$. 
With $\alpha$ of the form of \eqref{P2},
\[ Q_0=Q^2, \; \mbox{and} \; \frac{Q_0\Delta}{z} \leq \delta \leq Q_0, \]
we have, in a manner similar to Section 9 of \cite{Ba1}, for some absolute constant $c$, 
\begin{equation} \label{doh2}
 P \left( \alpha \right) \ll 1 + \delta^{-1} \int_{Q_0}^{4Q_0} \Omega( \delta, y) \dif y \ll 1 + \delta^{-1} \int_0^{\infty} e^{-y/Q_0} \Omega( \delta, y) \dif y,
\end{equation}
where
\[ \Omega ( \delta, y) = \sum\limits_{q\in \mathbbm{Z}} \ \phi\left(\frac{q-\sqrt{y}}{c\delta/\sqrt{Q_0}} \right) \sum\limits_{\substack{ m\in \mathbbm{Z} \\ m\equiv -bq^2 \pmod{r}}} \phi\left(\frac{m-yrz}{8\delta rz}\right) {\rm d}y. \]
Applying Poisson summation, Lemma~\ref{poisum}, twice in the same way as in Section 9 of \cite{Ba1}, we have
\begin{equation} \label{doh3}
 P(\alpha) \ll \frac{\delta z}{\sqrt{Q_0}} \left| \sum\limits_{j\in \mathbbm{Z}} 
\frac{\hat{\phi}(8j\delta z)}{r^*} \sum\limits_{l\in \mathbbm{Z}}
\hat{\phi}\left(\frac{cl\delta}{r^*\sqrt{Q_0}}\right) 
G(j^*b, l ; r^*) E(j,l) \right|,
\end{equation}
where 
\[ r^*=\frac{r}{(j,r)}, \; j^*=\frac{j}{(j,r)}, \]
$G(j^*b, l ; r^*)$ is defined as in \eqref{gausssumdef} and
\[ E(j,l) = \int_0^{\infty} e^{-y/Q_0} e\left(jyz-l\cdot\frac{\sqrt{y}}{r^*}\right) \ {\rm d}y. \]

\section{Evaluation of the exponential integrals}

We shall prove the following lemma in this section.
\begin{lemma}
Suppose that $N \geq 1$, $A, B \in \rear$. Then we have
\begin{enumerate}

\item If $B \neq 0$ and  $A/B > 0$, then 
\begin{equation} \label{intweightest1}
 \int_0^{\infty} e^{-y/N} e(Ay-B\sqrt{y}) \dif y = e \left( \frac{1}{8}
						     \right)
 \frac{2B}{(2A)^{3/2}} \exp \left( - \frac{B^2}{4A^2N} \right) e \left(
			     -\frac{B^2}{4A} \right) + O \left( \frac{1}{|A|}+\frac{1}{\sqrt{|A|}|B|} \right);
\end{equation}

\item If $B \neq 0$ and $A/B \leq 0$, then
\begin{equation} \label{intweightest2}
 \int_0^{\infty} e^{-y/N} e(Ay-B\sqrt{y}) \dif y \ll \frac{\sqrt{N}}{|B|} ;
\end{equation}

\item If $B=0$ and $A\neq 0$, then
\begin{equation} \label{intweightest3}
\int_0^{\infty} e^{-y/N} e(Ay-B\sqrt{y}) \dif y \ll \frac{1}{|A|}.
\end{equation}

\end{enumerate}
\end{lemma}

\begin{proof}
In (1) of the lemma, we may assume, without the loss of generality, that $B> 0$ and $A > 0$.  Then the left-hand side of \eqref{intweightest1} is, after the change of variables $x=\sqrt{y}$,
\begin{equation} \label{newint} 
2 \int_0^{\infty} e^{-x^2/N} e (Ax^2-Bx) x \dif x. 
\end{equation}
We split the above integral into 
\[ \int_0^a + \int_a^b + \int_b^{\infty}, \]
where $0<a<B/(4A)<B/A<b$. We observe that 
\begin{equation}\label{miniintegral} 
\lim\limits_{a\rightarrow 0} \int_0^a = 0, \; \mbox{and} \; \lim\limits_{b\rightarrow \infty} \int_b^{\infty} = 0.
\end{equation}
By the lemma concerning stationary phase, Lemma~\ref{Ivictheorem}, with the choices
\[ f(x)= F(x)=Ax^2, \; k=-B, \; g(x)=xe^{-\frac{x^2}{N}}, \; \mu(x)=\frac{x}{2}, \; \mbox{and} \; G(x)=x, \]
we have
\begin{equation} \label{hauptintegral}
 \int_a^b = e \left( \frac{1}{8} \right)  \frac{B}{(2A)^{3/2}} \exp \left( - \frac{B^2}{4A^2N} \right) e \left( -\frac{B^2}{4A} \right) + O \left( \frac{1}{A}+ \frac{1}{\sqrt{A}B}+\frac{b}{|2Ab-B|}+\frac{a}{|2Aa-B|}\right).
\end{equation}
Now \eqref{miniintegral} and \eqref{hauptintegral} imply (1) of the lemma upon taking $a$ and $b$ to zero and infinity, respectively.\newline

To prove (2) of the lemma, we split the integral in (\ref{newint}) into
\[ \int_0^b + \int_b^{\infty} \]
with large $b$,
apply Lemma~\ref{gkest} to the first integral after breaking up the 
interval $[0,b]$ into two subintervals on which $x\exp(-x^2/N)/(2Ax-B)$
is monotonic and using the second identity in \eqref{miniintegral}. \newline

Finally, if $B=0$, then the integral on the left-hand side of (\ref{intweightest3}) is
\[ \int_0^{\infty} \exp \left( y \left(2 \pi i A-\frac{1}{N}\right) \right) \dif y = \left( 2\pi iA -\frac{1}{N}\right)^{-1} \ll \frac{1}{|A|}. \]
Hence we have proved the lemma.
\end{proof}

\section{Treatment of the simple cases}

We now estimate the contributions of the simple parts of \eqref{doh3}.  The contribution of \eqref{doh3} with $j=l=0$ is, taking note that $r^*=1$ if $j=0$,
\begin{equation} \label{j=l=0}
\ll \delta z \sqrt{Q_0}.
\end{equation}
The contribution of $j \neq 0$ and $l=0$ is
\begin{equation} \label{jneq0l=0}
\ll \frac{\delta}{\sqrt{Q_0}} \sum\limits_{1\le j \le 1/(8\delta z)} \frac{1}{j\sqrt{r^*}} \ll \frac{\delta}{\sqrt{Q_0}} (Q_0 \Delta^{-1})^{\varepsilon} r^{-1/2},
\end{equation}
where the first of the above inequalities comes by the virtue of \eqref{intweightest3} and Lemma~\ref{gausssumest}, and the second arrives by the following estimate.
\begin{equation*}
\sum\limits_{1\le j \le 1/(8\delta z)} \frac{1}{j\sqrt{r^*}}= \frac{1}{\sqrt{r}}
\sum\limits_{t\vert r} \sqrt{t} \sum_{\substack{1\le j \le 1/(8\delta z) \\ (r,j)=t}} \frac{1}{j} \ll \frac{(Q_0 \Delta^{-1})^{\varepsilon}}{\sqrt{r}} 
\sum\limits_{t| r} \frac{1}{\sqrt{t}} \ll (Q_0 \Delta^{-1})^{2\varepsilon} r^{-1/2} .
\end{equation*}

If $l\not=0$ and $j/l \leq 0$, then it suffices to only consider that case in which $j \geq 0$ since the other case is similar and satisfies the same upper bound.  The contribution in question is,
by the virtue of \eqref{intweightest2} and Lemma~\ref{gausssumest},
\[  \ll \delta z\sum_{0 \leq j \leq (8\delta z)^{-1}} \sum_{\frac{-r^* \sqrt{Q_0}}{2c\delta} < l <0} \frac{\sqrt{r^*}}{|l|} \ll \delta z (Q_0 \Delta^{-1})^{\varepsilon} \sum_{0 \leq j \leq (8\delta z)^{-1}} \sqrt{r^*} \ll (Q_0 \Delta^{-1})^{\varepsilon} (\delta z + \sqrt{r}), \]
where we have used $r^*=1$ if $j=0$.  
Consequently, the total contribution to \eqref{doh3} of the above cases is 
\begin{equation} \label{simplecaseest}
 \ll (Q_0 \Delta^{-1})^{\varepsilon} \left( \sqrt{r} + \delta z \sqrt{Q_0} + \frac{\delta}{\sqrt{Q_0 r}} \right).
\end{equation}

\section{Transformation of the remaining terms}

We now consider the critical case when $j/l>0$. Then it again suffices to only consider that case in which $j>0$ since the other case is similar and satisfies the same upper bound.  By the virtue of \eqref{intweightest1}, Lemma~\ref{gausssumtrans} and Lemma~\ref{gausssumest}, the contribution in question is majorized by the sum of
\begin{equation} \label{j>0,l>01}
 \frac{\delta z}{\sqrt{Q_0}} \left| \sum\limits_{j>0} 
\sum\limits_{l>0}
g(j,l) G(j^*b,\chi(l),r^*) e \left( \frac{\overline{j^*b}(l^2-\chi(l))}{4r^*}+\frac{l^2}{4jz(r^*)^2}\right)
 \right|,
\end{equation}
and
\begin{equation} \label{j>0,l>02}
\frac{\delta z}{\sqrt{Q_0}} \sum\limits_{0<j<1/(8\delta z)}  
\sum\limits_{0<l<\frac{r^*\sqrt{Q_0}}{c\delta}} \frac{1}{\sqrt{r^*}}\left(
\frac{1}{jz}+\frac{r^*}{\sqrt{jz}l}\right),
\end{equation}
where 
\[ g(j,l) = \frac{\hat{\phi}(8j\delta z)}{r^*} \hat{\phi}\left(\frac{cl\delta}{r^*\sqrt{Q_0}}\right) \frac{l}{(jz)^{3/2}r^*} \exp \left( \frac{-l^2}{4(jz)^2Q_0(r^*)^2} \right), \]
and
\[ \chi(l) = \left\{ \begin{array}{ll} 0, & \mbox{if} \; 2|l; \\  1,  
& \mbox{if} \; 2 \nmid l. \end{array} \right. \]
Since $\delta\le Q_0$ and $r^* \leq r$, we find that \eqref{j>0,l>02} is 
majorized by
\[ \ll (Q_0\Delta^{-1})^{\varepsilon}\sqrt{r}. \] 

Now we break up the inner-most sum in (\ref{j>0,l>01}) into a sum over the even 
$l$'s and a sum over the odd $l$'s. In the following, we deal only with the contribution of the even $l$'s.  The contribution of the odd $l$'s can be estimated in a similar way and satisfies the same upper bound.  Using Lemma~\ref{gausssumest}, the part of \eqref{j>0,l>01} with $2|l$ is
\begin{equation} \label{beforeintpart}
 \ll \frac{\delta}{\sqrt{Q_0z}} \sum_{0 < j < (8\delta z)^{-1}} \frac{1}{(jr^*)^{3/2}} \left| \sum_{0<l< \frac{r^*\sqrt{Q_0}}{2c\delta}} \Phi (j, 2l)  e \left( \frac{\overline{j^*b}l^2}{r^*}+\frac{l^2}{jz(r^*)^2} \right) \right|,
\end{equation}
where
\[ \Phi (j, 2l) = \left( 1- \frac{2cl \delta}{r^* \sqrt{Q_0}} \right) l \exp \left( \frac{-l^2}{Q_0 (jzr^*)^2} \right).  \]
We now apply partial summation to the inner-most sum of \eqref{beforeintpart}.  It becomes
\begin{equation} \label{afterintpart}
-\int_0^{\frac{r^*\sqrt{Q_0}}{2c\delta}} \left(\frac{{\rm d}}{{\rm d}x}
\Phi(j,2x)\right) \sum_{0<l<x} e \left( \frac{\overline{j^*b}l^2}{r^*}+\frac{l^2}{jz(r^*)^2} \right) \dif x.
\end{equation}
Note that the boundary terms vanish in the partial summation.  Now we set
\[ D(j) = \min \left( \sqrt{Q_0} jzr^* (Q_0\Delta^{-1})^{\varepsilon}, \frac{r^*\sqrt{Q_0}}{2c\delta} \right) . \]
Then we break up the integral in \eqref{afterintpart} into integrals over the intervals $[0,D(j)]$ and $[D(j), (r^*\sqrt{Q_0})/(2c\delta) ]$.  The former of the two integrals is
\begin{equation} \label{firstint}
\ll \max_{0 < L \leq D(j)} \left| \sum_{0 < l \leq L} e \left( \frac{\overline{j^*b}l^2}{r^*}+\frac{l^2}{jz(r^*)^2} \right) \right| \int_0^{D(j)} \left| 
\frac{{\rm d}}{{\rm d}x} \Phi(j,2x) \right| \dif x,
\end{equation}
and the integral in \eqref{firstint} is
\[ \ll \max_{0 \leq x \leq D(j)} \Phi(j, 2x) \leq D(j) \leq \sqrt{Q_0} jzr^* (Q_0\Delta^{-1})^{\varepsilon}; \]
and the latter is
\begin{equation} \label{secondint}
\ll \max_{D(j) \leq L \leq \frac{r^*\sqrt{Q_0}}{2c\delta}} \left| \sum_{0 < l \leq L} e \left( \frac{\overline{j^*b}l^2}{r^*}+\frac{l^2}{jz(r^*)^2} \right) \right| \int_{D(j)}^{\frac{r^*\sqrt{Q_0}}{2c\delta}} \left| 
\frac{{\rm d}}{{\rm d}x} \Phi(j,2x) \right| \dif x,
\end{equation}
and the integral in \eqref{secondint} is, for any $C>0$,
\[ \ll (Q_0 \Delta^{-1})^{-C} .\]
Therefore, \eqref{afterintpart} is, for any $C>0$,
\begin{equation} \label{secondinpart}
\ll \max_{0 < L \leq D(j)} \left| \sum_{0 < l \leq L} e \left( \frac{\overline{j^*b}l^2}{r^*}+\frac{l^2}{jz(r^*)^2} \right) \right| \sqrt{Q_0} jzr^* (Q_0\Delta^{-1})^{\varepsilon} + (Q_0 \Delta^{-1})^{-C}.
\end{equation}

Now we set
\begin{equation} \label{vonbismarck}
 k^* = \frac{1}{zr^*r}.
\end{equation}
$k^*$ is a positive integer by \eqref{otto} and the fact that $r^*|r$.  We further assume that
\[ \overline{b} \equiv -a \pmod{r^*}, \; \mbox{with} \; 1 \leq a \leq r^*. \]
We write
\[ \frac{\overline{j^*}}{r^*} \equiv - \frac{\overline{r^*}}{j^*} + \frac{1}{j^*r^*} \pmod{1}. \]
Then
\[ e \left( \frac{\overline{j^*b}l^2}{r^*}+\frac{l^2}{jz(r^*)^2} \right) = e \left( \frac{-al^2}{j^*r^*} \right) e \left( \frac{a\overline{r^*}+k^*}{j^*}l^2 \right), \]
where we have used the relation $r^*j=rj^*$. We will use this relation
in several places of the remainder of this paper. 
We now apply partial summation to the sum over $l$ in \eqref{secondinpart} with the above inserted into the appropriate place.  We get that the sum in question is
\[ \ll \left( 1+\frac{L^2}{j^*} \right) \max_{0<x \leq L} \left| \sum_{0 <l \leq x} e \left( \frac{a\overline{r^*}+k^*}{j^*}l^2 \right) \right|. \]
Now combining everything thus far, we have that the contribution from $j$'s and $l$'s with $j/l >0$ is
\begin{equation} \label{afteralltrans}
\ll (Q_0 \Delta^{-1})^{\varepsilon} \sqrt{r} + (Q_0 \Delta^{-1})^{\varepsilon} \delta \sqrt{z} \sum_{0<j<(8\delta z)^{-1}} \frac{1}{\sqrt{r^*j}} \max_{0 < x \leq D(j)} \left| \sum_{0<l\leq x} e \left( \frac{a\overline{r^*}+k^*}{j^*}l^2 \right) \right|,
\end{equation}
where we have used the estimates
\[ \frac{L^2}{j^*} \leq \frac{D(j)^2}{j^*} \leq \frac{Q_0(jzr^*)^2}{j^*}(Q_0 \Delta^{-1})^{2\varepsilon} \] 
and
\[ \frac{Q_0(jzr^*)^2}{j^*} = Q_0(jz)^2 r^* \frac{r}{j} \leq Q_0jz^2r^2 \leq Q_0j\Delta \leq \frac{Q_0\Delta}{8\delta z} \leq 1. \]

\section{Application of Weyl shift}

Now the inner-most sum in \eqref{afteralltrans} can be estimated again using Weyl shift, Lemma~\ref{weylk}.  We have
\[ \left| \sum_{0<l\leq x} e \left( \frac{a\overline{r^*}+k^*}{j^*}l^2 \right) \right|^2 \ll x + \sum_{1\leq l \leq x} \min \left( x, \left\| \frac{a\overline{r^*}+k^*}{j^*}2l\right\|^{-1} \right). \]
Therefore, we infer that the double sum in the last term in \eqref{afteralltrans} is,
\begin{equation} \label{beforecauchy}
\begin{split}
\ll \frac{1}{\sqrt{r}}  \sum_{t|r} &\sum_{\substack{0<j^*<(8\delta z t)^{-1} \\ (r^*,j^*)=1}} \sqrt{\frac{D(j^*t)}{j^*}} \\
&+ \frac{1}{\sqrt{r}}  \sum_{t|r} \sum_{\substack{0<j^*<(8\delta z t)^{-1} \\ (r^*,j^*)=1}} \frac{1}{\sqrt{j^*}} \left( \sum_{1\leq l \leq 2D(j^*t)} \min \left( D(j^*t), \left\| \frac{a\overline{r^*}+k^*}{j^*}l \right\|^{-1} \right) \right)^{\frac{1}{2}},
\end{split}
\end{equation}
where $r^*=r/t$ and $j^*=j/t$.  The first term in \eqref{beforecauchy} is
\[ \ll (Q_0\Delta^{-1})^{2\varepsilon} \frac{Q_0^{1/4}}{\delta \sqrt{z}}, \]
where we have used the estimate
\begin{equation} \label{djstar}
 D(j^*t) \leq \sqrt{Q_0} jzr^* (Q_0\Delta^{-1})^{\varepsilon} = \sqrt{Q_0} j^*zr (Q_0\Delta^{-1})^{\varepsilon}. 
\end{equation}
Now we apply Cauchy's inequality to the inner double sum of the second term of \eqref{beforecauchy}.  We obtain that its square is
\begin{equation} \label{aftercauchy}
 \ll (Q_0\Delta^{-1})^{\varepsilon} \sum_{\substack{0<j^*<(8\delta z t)^{-1} \\ (r^*,j^*)=1}} \sum_{1\leq l \leq 2D(j^*t)} \min \left( D(j^*t), \left\| \frac{a\overline{r^*}+k^*}{j^*}l \right\|^{-1} \right).
\end{equation}
Now to estimate \eqref{aftercauchy}, it suffices to estimate the number of solutions $(j^*,l, h)$ with
\begin{equation} \label{boundjl} 
J < j^* \leq 2J, \; (r^*,j^*)=1, \; 1 \leq l \leq 4 (Q_0\Delta^{-1})^{\varepsilon} Q_0^{1/2}Jzr, \; \mbox{and} \; H \leq h \leq 2H
\end{equation}
to the equation
\[ \left\| \frac{(a\overline{r^*}+k^*)l}{j^*} \right\| = \frac{h}{j^*}. \]
The number in question is majorized by the number of solutions in $(j^*,l, h)$ of the congruence
\begin{equation*}
 \pm (k^*+a\overline{r^*}) l \equiv h \pmod{j^*}
\end{equation*}
with $j^*$, $l$ and $h$ satisfying \eqref{boundjl}.  The above congruence is equivalent to
\begin{equation} \label{congruenceeq}
 \pm (k^*r^*+a) l -hr^* \equiv 0 \pmod{j^*},
\end{equation}
since $(r^*,j^*)=1$.  First, if $ \pm (k^*r^*+a) l = hr^*$, then every $(j^*, l, h)$ with $J < j^* \leq 2J$ is a solution to \eqref{congruenceeq}.  Moreover the number of solutions to
\[ \pm (k^*r^*+a) l = hr^* \]
is at most $H/(k^*r^*+a)$, since $(k^*r^*+a, r^*)=1$.  Hence the total number of solutions in this case is
\begin{equation*}
\ll J\frac{H}{k^*r^*+a} \leq JH zr
\end{equation*}
since by \eqref{vonbismarck}
\[ k^*r^*+a > k^*r^* = \frac{1}{zr} . \]
Second if $ \pm (k^*r^*+a) l \neq hr^*$, then we fix $l$ and $h$.  Hence the number of $j^*$'s such that $(j^*,l,h)$ is a solution is at most
\[ \tau(|(k^*r^*+a) l - hr^*|) + \tau ((k^*r^*+a) l + hr^*) \ll (Q_0\Delta^{-1})^{\varepsilon} . \]
Therefore, the total number of solutions in this case is
\[ \ll Q_0^{1/2}Jzr H(Q_0\Delta^{-1})^{\varepsilon}. \]
Now combining the two cases, we get that the number of solutions $(j^*,l,h)$ for \eqref{congruenceeq} satisfying \eqref{boundjl} is
\begin{equation} \label{saladin}
\ll Q_0^{1/2}Jzr H(Q_0\Delta^{-1})^{\varepsilon}.
\end{equation}

We now write 
\[ \sum\limits_{j^*,l,h}= \mathop{\sum_{J < j^* \leq 2J} \ \sum_{0<l\leq 2D(j^*t)} \ \sum_{H \leq h \leq 2H}}_{\substack{\|(a\overline{r^*}+k^*)l/j^*\| = h/j^* \\ (j^*,r^*)=1}} \ \ \ \ \ \ \mbox{and} \ \ \ \ \ \ \sum\limits_{j^*,l}=
\mathop{\sum_{J < j^* \leq 2J} \ \sum_{0<l\leq 2D(j^*t)}}_{j^*|(a\overline{r^*}+k^*)l, \;(j^*,r^*)=1}. \]
To estimate \eqref{aftercauchy}, it suffices to estimate the following.
\begin{eqnarray}
\nonumber & & \mathop{\sum_{J < j^* \leq 2J} \ \sum_{0<l\leq 2D(j^*t)} \ \sum_{0 \leq h \leq 2J}}_{\substack{\|(a\overline{r^*}+k^*)l/j^*\| = h/j^* \\ (j^*,r^*)=1}} \min \left( D(j^*t), \frac{j^*}{h} \right) \\
\nonumber & &\ll  (Q_0\Delta^{-1})^{\varepsilon} \max_{1 \leq H \leq J} 
\sum_{j^*,l,h} \min \left( D(j^*t), \frac{j^*}{H} \right) + 
(Q_0\Delta^{-1})^{\varepsilon}  \sum_{j^*,l} D(j^*t) \\
\label{lenin}& &  \ll (Q_0\Delta^{-1})^{2\varepsilon} \left( \max_{1 \leq H \leq (\sqrt{Q_0}zr)^{-1}} \sum_{j^*,l,h} \sqrt{Q_0}Jzr + \max_{(\sqrt{Q_0}zr)^{-1}<H \leq J} \sum_{j^*,l,h} \frac{J}{H} + Q_0(Jzr)^2 \right)
\end{eqnarray}
where we have used the estimate in \eqref{djstar}.  Now using \eqref{saladin}, \eqref{lenin} is
\[  \ll (Q_0\Delta^{-1})^{\varepsilon} \left( J^2 \sqrt{Q_0}zr + Q_0(Jzr)^2 \right). \]
Now combining everything, we get that \eqref{afteralltrans} is
\begin{equation} \label{marx}
 \ll (Q_0\Delta^{-1})^{\varepsilon} \left( Q_0^{1/4} + \sqrt{Q_0} \Delta^{1/4}+ \sqrt{r} \right),
\end{equation}
where we have used \eqref{P2}.

\section{A second estimate for $P(\alpha)$}

We first assume that $z\ge \Delta$.  
Combining \eqref{simplecaseest} and \eqref{marx} and using \eqref{doh3} with 
\[ \delta = \frac{Q_0 \Delta}{z}, \]
we get that
\[ P ( \alpha) \ll (Q_0\Delta^{-1})^{\varepsilon} \left( Q_0^{1/4} + Q_0^{1/2} \Delta^{1/4} + Q_0^{3/2} \Delta + \sqrt{r} + \frac{Q_0^{1/2} \Delta}{r^{1/2} z} \right). \]
On the other hand, by the first inequality of Lemma 6 of \cite{Ba1}, we have
\[ P\left( \alpha \right) \ll \Delta^{-\varepsilon} \left(1+Q_0rz+Q_0^{3/2}\Delta\right). \]
Combining the above two inequalities, we get
\[  P ( \alpha) \ll \left( \frac{Q_0}{\Delta} \right)^{\varepsilon} \left( Q_0^{1/4} + Q_0^{1/2} \Delta^{1/4} + Q_0^{3/2} \Delta + \sqrt{r} + \min \left( \frac{\sqrt{Q_0} \Delta}{\sqrt{r} z}, Q_0rz \right) \right). \]
Now as in Section 7 of \cite{Ba1}, we have
\[ \min \left( \frac{Q_0^{1/2} \Delta}{r^{1/2} z}, Q_0rz \right) \leq Q_0^{3/4} \Delta^{1/2}r^{1/4}. \]
Recalling that $Q_0=Q^2$, and using (\ref{verysmalldelta}) in the case when
$\Delta/2\le z\le \Delta$,
we have proved the following.

\begin{theorem}
Let $\varepsilon >0$ be given and $\alpha$ satisfy \eqref{alphacond}.  We have
\begin{equation} \label{secondestPalpha}
P(\alpha) \ll \left( \frac{Q}{\Delta} \right)^{\varepsilon} \left( Q^{1/2} + Q \Delta^{1/4} + Q^3 \Delta + \sqrt{r} + Q^{3/2} \Delta^{1/2}r^{1/4} +
Q^2\Delta r\right),
\end{equation}
where the implied constant depends on $\varepsilon$ alone.
\end{theorem}

\section{Proof of Theorem~\ref{squaremod1}}

Finally, we are ready to prove Theorem~\ref{squaremod1}.  Noting $r \leq \Delta^{-1/2}$ by \eqref{P2}, we use \eqref{secondestPalpha} if $r \leq Q$ and \eqref{firstestPalpha} if $r >Q$ to get the following.
\[ P(\alpha) \ll \left( \frac{Q}{\Delta} \right)^{\varepsilon} \left( Q^3\Delta + Q^{7/4} \Delta^{1/2} + Q \Delta^{1/4} +Q^{1/2}\right). \]
Choosing $\Delta = N^{-1}$ and using Lemma~\ref{modified} and Lemma~\ref{engels} after dividing the outer-most sum in \eqref{aim} 
into sums over
dyadic intervals, we get that \eqref{aim} is bounded by
\[ \ll \left( QN\right)^{\varepsilon} \left( Q^3 +  Q^{7/4} N^{1/2} + Q N^{3/4} + Q^{1/2}N \right)Z. \]
If $Q\ge N^{2/5}$ then
\[ Q^{7/4} N^{1/2} \leq Q^3, \; \mbox{and} \; Q^{1/2}N \leq Q^3. \]
On the other hand, if $Q< N^{2/5}$ then
\[ Q^{7/4} N^{1/2} \leq Q^{1/2}N . \]
Moreover, if $Q \leq N^{1/2}$, then 
\[ QN^{3/4} \leq Q^{1/2}N. \]
Therefore, if $Q \leq N^{1/2}$, we have that the left-hand side of \eqref{squaremod1eq} is
\[ \ll \left( QN\right)^{\varepsilon}  \left( Q^3 + Q^{1/2}N \right) Z. \]
If $Q > N^{1/2}$ then the above majorizes the first quantity in \eqref{A4}.  Now comparing the above with \eqref{Baiertheorem}, we have the theorem.\\

{\bf Acknowledgement.}
This paper was written when the first and second-named authors held postdoctoral fellowships at the Department of Mathematics and Statistics at Queen's University and the Department of Mathematics at the University of Toronto, respectively.  The authors wish to thank these institutions for their financial support.

\bibliography{biblio}
\bibliographystyle{amsxport}

\vspace*{.7cm}

\noindent Department of Mathematics and Statistics, Queen's University \newline
University Ave, Kingston, ON K7L 3N6 Canada \newline
Email: {\tt sbaier@mast.queensu.ca} \newline

\noindent Department of Mathematics, University of Toronto \newline
40 Saint George Street, Toronto, ON M5S 2E4 Canada \newline
Email: {\tt lzhao@math.toronto.edu}
\end{document}